\documentclass[11pt]{article}

 \usepackage{setspace}
\usepackage{natbib}

\usepackage{latexsym}
\usepackage{calc}

\usepackage{amssymb}
\usepackage{amsmath}
\usepackage{graphicx}

\usepackage{psfrag}

\usepackage{float}

\usepackage{color}

\usepackage{nicefrac}

\usepackage{ntheorem}

\usepackage{hyperref}

\usepackage{tikz}
\usetikzlibrary{backgrounds}

\newcounter{hours}\newcounter{minutes}

\def\nr{\par \noindent}

\def\beq{\begin{equation}}
\def\eeq{\end{equation}}

\newtheorem{theorem}{Theorem}

\newtheorem{lemma}{Lemma}
\newtheorem{corollary}{Corollary}

\newtheorem{proposition}{Proposition}
\newtheorem{assumption}{Assumption}
\newtheorem{definition}{Definition}

\newtheorem{example}{Example}
\newtheorem{remark}{Remark}
\newcommand{\proof}{\bf Proof: \rm \nr}

\newcommand{\qed}{\hfill $\Box$ \nr \medskip}

\def\ba{\begin{array}}
\def\ea{\end{array}}
\def\beann{\begin{eqnarray*}}
\def\eeann{\end{eqnarray*}}
\def\bea{\begin{eqnarray}}
\def\eea{\end{eqnarray}}

\def\BT{\begin{theorem}}
\def\ET{\end{theorem}}
\def\BL{\begin{lemma}}
\def\EL{\end{lemma}}
\def\BC{\begin{corollary}}
\def\EC{\end{corollary}}
\def\BE{\begin{example}}
\def\EE{\end{example}}
\def\BD{\begin{definition}}
\def\ED{\end{definition}}
\def\BR{\begin{remark}}
\def\ER{\end{remark}}
\def\BAS{\begin{assumption}}
\def\EAS{\end{assumption}}
\def\BI{\begin{itemize}}
\def\EI{\end{itemize}}

\def\BMP{\begin{minipage}{9.5cm}}
\def\EMP{\end{minipage}}
\def\MPT{\begin{minipage}{11.5cm}}
\def\EPT{\end{minipage}}

\def\R{\mathbb{R}}

\newcommand*\samethanks[1][\value{footnote}]{\footnotemark[#1]}
\newcommand{\NNorm}[2]{\left\Vert {#1} \right\Vert_{#2}}
\newcommand{\AV}[1]{\left\vert {#1} \right\vert}

\title{On nondegenerate M-stationary points for \\ sparsity constrained nonlinear optimization}
\author{
S. L\"ammel
\thanks{
Department of Mathematics, Chemnitz University of Technology,
Reichenhainer Str. 41, 09126
Chemnitz, Germany; e-mail: sebastian.laemmel@mathematik.tu-chemnitz.de, vladimir.shikhman@mathematik.tu-chemnitz.de (corresponding author).
 } \and V. Shikhman\samethanks[1]
}

\begin{document}
\maketitle
\vspace{-5ex}
\abstract{ 
We study sparsity constrained nonlinear optimization (SCNO) from a
topological point of view. Special focus will be on M-stationary points from \cite{burdakov:2016}. We introduce nondegenerate M-stationary
points and define their M-index. We show that all M-stationary points are generically nondegenerate. In particular, the sparsity constraint is active at all local minimizers of a generic SCNO. Some relations to other stationarity concepts, such as S-stationarity, basic feasibility, and CW-minimality, are discussed in detail. By doing so, the issues of instability and degeneracy of  points due to different stationarity concepts are highlighted.
The concept of M-stationarity allows to adequately describe the global structure of SCNO along the lines of Morse theory. For that, we study topological changes of lower level sets while passing an M-stationary point. As novelty for SCNO, multiple cells of dimension equal to the M-index are needed to be attached. This intriguing fact is in strong contrast with other optimization problems considered before, where just one cell suffices. As a consequence, we derive a Morse relation for SCNO, which relates the numbers of local minimizers and M-stationary points of M-index equal to one. The appearance of such saddle points cannot be thus neglected from the perspective of global optimization. Due to the multiplicity phenomenon in cell-attachment, a saddle point may lead to more than two different local minimizers. We conclude that the relatively involved structure of saddle points is the source of well-known difficulty if solving SCNO to global optimality. 
}

\vspace{2ex}

\section{Introduction}
\label{sec:intro}
We consider the sparsity constrained nonlinear optimization:
\[
\label{eq:SCNO}
\mbox{SCNO:}\quad \min_{x \in \R^n}\,\, f(x)\quad \mbox{s.\,t.} \quad \NNorm{x}{0} \le s,
\]
where the so-called $\ell_0$ "norm" counts non-zero entries of $x$:
\[
\NNorm{x}{0}= \AV{\left\{i\in \{1,\ldots,n\} \,\left\vert\, x_i \ne 0 \right.\right\}},
\]
the objective function $f \in C^2\left(\R^n,\R\right)$ is twice continuously differentiable, and $s \in\{0,1, \ldots,n-1\}$ is an integer.
The difficulty of solving SCNO comes from the combinatorial nature of the sparsity constraint $\NNorm{x}{0} \le s$. The requirement of sparsity is however motivated by various applications, such as compressed sensing, model selection, image processing etc. We refer e.\,g. to \cite{donoho:2006}, \cite{tibshirani:1996}, and \cite{shechtman:2011} for further details on the relevant applications.

In the seminal paper \cite{beck:2013}, necessary optimality conditions for SCNO have been stated. Namely, the notions of basic feasibility (BF-vector), $L$-stationarity and CW-minimality have been introduced and studied there. 
Note that the formulation of $L$-stationarity mimics the techniques from convex optimization by using the orthogonal projection on the SCNO feasible set. The notion of CW-minimum incorporates the coordinate-wise optimality along the axes. Based on both stationarity concepts, algorithms that find points satisfying these conditions have been developed. Those are the iterative hard thresholding method, as well as the greedy and partial sparse-simplex methods. In a series of subsequent papers \cite{beck:2016, beck:2018} elaborated the algorithmic approach for SCNO which is based on $L$-stationarity and CW-minimality. 

Another line of research started with \cite{burdakov:2016}, where additionally smooth equality and inequality constraints have been incorporated into SCNO. For that, the authors coin the new term of mathematical programs with cardinality constraints (MPCC). The key idea in \cite{burdakov:2016} is to provide a mixed-integer formulation whose standard relaxation still has the same solutions as MPCC. For the relaxation the notion of S-stationary points is proposed. S-stationarity corresponds to the standard Karush-Kuhn-Tucker condition for the relaxed program. The techniques applied follow mainly those for mathematical programs with complementarity constraints.
In particular, an appropriate regularization method for solving MPCC is suggested. The latter is proved to converge towards so-called M-stationary points. M-stationarity corresponds to the standard Karush-Kuhn-Tucker condition of the tightened program, where zero entries of an MPCC feasible point remain locally vanishing. Further research in this direction is presented in a series of subsequent papers \cite{cervinka:2016}, \cite{bucher:2018}.

The goal of this paper is the study of SCNO from a topological point of view. The topological approach to optimization has been pioneered by \cite{jongen:1977}, \cite{Jongen:2000} for nonlinear programming problems, and successfully developed for mathematical programs with complementarity constraints, mathematical problems with vanishing constraints, general semi-infinite programming, bilevel optimization, semi-definite programming, disjunctive programming etc., see e.\,g. \cite{Shikhman:2012} and references therein. The main idea of the topological approach is to identify stationary points which roughly speaking induce the global structure of the underlying optimization problem. The stationary points include minimizers, but also all kinds of saddle points -- just in analogy to the unconstrained case.   
It turns out that for SCNO the concept of M-stationarity from \cite{burdakov:2016} is the adequate stationarity concept at least from  the topological perspective. We outline our main findings and results:

\begin{itemize}
    \item[1.] We introduce nondegenerate M-stationary points along with their associated M-indices. The latter subsume as usual the quadratic part -- the number of negative eigenvalues of the objective's Hessian restricted to non-vanishing variables. As novelty, the sparsity constraint provides an addition to the M-index, namely, the difference between the bound and the current number of non-zero variables at a nondegenerate M-stationary point. We prove that all M-stationary points are generically nondegenerate. In particular, it follows that all local minimizers of SCNO are nondegenerate with vanishing M-index, hence, the sparsity constraint is active. Note that M-stationary points with non-vanishing M-index correspond to saddle points. The local structure of SCNO around a nondegenerate M-stationary point is fully described just by its M-index, at least up to a differentiable change of coordinates. 
    \item[2.] We thoroughly discuss the relation of M-stationarity to S-stationarity, basic feasibility, and CW-minimality for SCNO. It turns out that nondegenerate M-stationary points may cause degeneracies of S-stationary points viewed as Karush-Kuhn-Tucker-points for the relaxed problem. Moreover, even under the cardinality constrained second-order sufficient optimality condition from \cite{bucher:2018} assumed to hold at an S-stationary point, the corresponding M-stationary point does not need to be a nondegenerate local minimizer for SCNO.
    As for CW-minima, we show that they are not stable with respect to data perturbations in SCNO. After an arbitrarily small $C^2$-perturbation of $f$ a locally unique CW-minimum may bifurcate into multiple CW-minima. More importantly, this bifurcation unavoidably causes the emergence of M-stationary points, being different from the CW-minima. Despite of this instability phenomenon, if a BF-vector and, hence, CW-minimum, happens to be nondegenerate as an M-stationary point, then the sparsity constraint is necessarily active. 
    \item[3.] We use the concept of M-stationarity in order to describe the global structure of SCNO. To this aim the study of topological properties of its lower level sets is undertaken. As in the standard Morse theory, see e.\,g. \cite{milnor:1963}, \cite{goresky:1988}, we focus on the topological changes of the lower level sets as their levels vary. Appropriate versions of deformation and cell-attachment theorems are shown to hold for SCNO. Whereas the deformation is standard, the cell-attachment reveals an essentially new phenomenon not observed in nonsmooth optimization before.  
    In SCNO, {\it multiple} cells of the
same dimension need to be attached, see Theorem \ref{thm:cell-a}. 
To determine the number of these attached cells turns out to constitute a challenging combinatorial problem from algebraic topology, see Lemma \ref{lem:cat}.  
\item[4.] As a consequence of proposed Morse theory, we derive a Morse relation for SCNO, which relates the numbers of local minimizers and M-stationary points of M-index equal to one. The appearance of such saddle points cannot be thus neglected from the perspective of global optimization. As novelty for SCNO, a saddle point may lead to more than two different local minimizers. This is in strong contrast with other nonsmooth optimization problems studied before, see e.\,g. \cite{Shikhman:2012}, where a saddle point leads to at most two of them. We conclude that the relatively involved structure of saddle points is the source of well-known difficulty if solving SCNO to global optimality. 
\end{itemize}

The paper is organized as follows. In Section \ref{sec:m-stat} we discuss the notion of M-stationarity for SCNO. Section \ref{sec:rel} is devoted to the relation of M-stationarity to other stationarity concepts from the literature. In Section \ref{sec:crit} the global structure of SCNO is described within the scope of Morse theory.

Our notation is standard. The cardinality of a finite set $S$ is denoted by $|S|$. The $n$-dimensional Euclidean space is denoted by $\R^n$ with the coordinate vectors $e_i$, $i=1, \ldots,n$. 
For $J \subset \{1, \ldots, n\}$ we denote by $\mbox{conv}\left(e_j, j \in J\right)$ the convex combination of the coordinate vectors $e_j, j \in J$. Given a twice continuously differentiable function $f: \R^n \rightarrow \R$, $\nabla f$ denotes its gradient, and $D^2 f$ stands for its Hessian.

\section{M-stationarity}
\label{sec:m-stat}

For $0 \leq k \leq n$ we use the notation
\[
  \R^{n,k}= \left\{ x \in \R^n\, \left\vert \, \NNorm{x}{0} \le k \right.\right\}.
\]
Using the latter, the feasible set of SCNO can be written as
\[
  \R^{n,s}= \left\{ x \in \R^n\, \left\vert \, \NNorm{x}{0} \le s \right.\right\}.
\]
For a feasible point $x\in \R^{n,s}$ we define the following complementary index sets:
\[
I_0(x) = \left\{i\in \{1,\ldots,n\} \,\left\vert\, x_i = 0 \right.\right\}, \quad
  I_1(x) = \left\{i\in \{1,\ldots,n\} \,\left\vert\, x_i \ne 0 \right.\right\}.
\]
Without loss of generality, we assume throughout the whole paper that
at the particular point of interest $\bar x \in  \R^{n,s}$ with $\NNorm{\bar x}{0}=k$ it holds:
\[
  I_0\left(\bar x\right) = \left\{1,\ldots,n-k\right\}, \quad
  I_1\left(\bar x\right) = \left\{n-k+1,\ldots,n\right\}.
\]
Using this convention, the following local description of SCNO feasible set can be deduced. Let $\bar x \in \R^{n,s}$ be a feasible point for SCNO with $\NNorm{\bar x}{0}=k$. Then, there exist neighborhoods $U_{\bar x}$ and $V_0$ of $\bar x$ and $0$, respectively, such that under the linear coordinate transformation $\Phi(x)=x-\bar x$ we have locally:
\begin{equation}
\label{eq:diff}
\Phi\left(\R^{n,s}\cap U_{\bar x}\right) = \left( \R^{n-k,s-k}\times \R^k \right) \cap V_0, \quad \Phi(\bar x)=0.
\end{equation}{}
 


\begin{definition}[M-stationarity, \cite{burdakov:2016}]
\label{def:m-stat}
A feasible point $\bar x \in \R^{n,s}$ is called M-stationary for SCNO if 
\[
\frac{\partial f}{\partial x_i} \left( \bar x\right) = 0  \mbox{ for all } i \in I_1\left(\bar x\right).
\]
\end{definition}
%
Obviously, a local minimizer of SCNO is 
an M-stationary point. 


%

\begin{definition}[Nondegenerate M-stationarity]
\label{def:nondeg}
An M-stationary point $\bar x \in \R^{n,s}$ with $\left\| \bar x\right\|_0 =k$ is called nondegenerate if the following conditions hold:
\begin{itemize}
\item[] ND1: if $k < s$ then $\displaystyle \frac{\partial f}{\partial x_i} \left( \bar x \right) \ne 0$ for all $\displaystyle  i \in I_0\left(\bar x\right)$,
\item[] ND2: the matrix $\displaystyle \left(\frac{\partial^2 f}{\partial x_i \partial x_j}\left(\bar x \right) \right)_{i,j \in I_1\left(\bar x\right)}$ is nonsingular.
\end{itemize}
Otherwise, we call $\bar x$ degenerate.
\end{definition}

\begin{definition}[M-Index]
\label{def:mindex}
Let $\bar x \in \R^{n,s}$ be a nondegenerate M-stationary point with $\NNorm{\bar x}{0}=k$.
The  number of negative eigenvalues of the matrix $\displaystyle \left(\frac{\partial^2 f}{\partial x_i \partial x_j}\left(\bar x \right) \right)_{i,j \in I_1\left(\bar x\right)}$ is called its quadratic index ($QI$). 
The number $s-k+QI$ is called the M-index of $\bar x$. 
\end{definition}

\begin{theorem}[Morse-Lemma for SCNO]
\label{thm:morse}
Suppose that $\bar x$ is a nondegenerate M-stationary point for SCNO with  $\NNorm{\bar x}{0}= k$ and  quadratic index $QI$. Then, there exist neighborhoods $U_{\bar x}$ and $V_0$ of $\bar x$ and $0$, respectively, and a local $C^1$-coordinate system $\Psi: U_{\bar x} \rightarrow V_0$ of $\R^n$ around $\bar x$ such that:

\begin{equation}
\label{eq:normal}
    f\circ \Psi^{-1}(y)= f(\bar x) + \sum\limits_{i=1}^{n-k}y_i + \sum\limits_{j=n-k+1}^{n}\pm y_j^2,  
\end{equation}
where $y \in \R^{n-k,s-k}\times \R^k$. Moreover, there are exactly $QI$ negative squares in (\ref{eq:normal}).
\end{theorem}

\proof
Without loss of generality, we may assume $f\left(\bar x\right)=0$.
By using $\Phi$ from (\ref{eq:diff}), we put $\bar f:= f\circ \Phi^{-1}$ on the set $\left(\R^{n-k,s-k}\times \R^k \right)\cap V_0$.
At the origin we have:
\begin{itemize}
    \item[(i)] if $k < s$ then $\displaystyle \frac{\partial \bar f}{\partial y_i} \ne 0$ for all $i =1,\ldots, n-k$,
    \item[(ii)] $\displaystyle \frac{\partial \bar f}{\partial y_i} = 0$  for all $i =n-k+1, \ldots, n$,
\item[(iii)] the matrix $\displaystyle \left(\frac{\partial^2 \bar f}{\partial y_i \partial y_j} \right)_{i,j =n-k+1, \ldots, n}$ is nonsingular.
\end{itemize}
We denote $\bar f$ by $f$ again. Under the following coordinate transformations the set
$\R^{n-k,s-k}\times \R^k$ will be equivariantly transformed in itself.
We put $y=\left(Y_{n-k},Y^k\right)$, where $Y_{n-k}=(y_1,\ldots,y_{n-k})$ and $Y^{k}=(y_{n-k+1},\ldots,y_n)$. It holds:
\[
f\left(Y_{n-k},Y^k\right)= \int_0^1 \frac{d}{dt} f\left(tY_{n-k},Y^{k}\right)dt+ f\left(0,Y^{k}\right)=  \sum_{i=1}^{n-k}y_id_i(y)+f\left(0,Y^{k}\right),
\]
where 
\[
 d_i(y) = \int_0^1 \frac{\partial f}{\partial y_i} \left(tY_{n-k},Y^{k}\right)dt, \quad i=1,\ldots,n-k.
\]
Note that $d_i\in C^1, i=1,\ldots,n-k$. Due to (ii)-(iii), we may apply the standard Morse lemma on the $C^2$-function $ f\left(0,Y^{k}\right)$ without affecting the coordinates $Y_{n-k}$, see e.\,g. \cite{Jongen:2000}. The corresponding coordinate transformation is of class $C^1$.
Denoting the transformed functions again by $f$ and $d_i$, we obtain
\[
f(y) = \sum_{i=1}^{n-k}y_id_i(y) + \sum\limits_{j=n-k+1}^{n}\pm y_j^2.
\]

In case $k=s$, we need to consider $f$ locally around the origin on the set 
\[
\R^{n-k,s-k}\times \R^k = \R^{n-k,0}\times \R^k = \{0\}^{n-k}\times \R^k.
\]
Hence, $y_i=0$ for $i=1,\ldots, n-k$, and we immediately obtain the representation (\ref{eq:normal}).

In case $k<s$, (i) provides that $\displaystyle d_i(0)=\frac{\partial f}{\partial y_i}(0) \not = 0$, $i=1, \ldots, n-k$.
Hence, we may take
\[
 y_i d_i(y), i=1,\ldots,n-k, \quad y_j,  j=n-k+1,\ldots, n
\]
as new local $C^1$-coordinates by a straightforward application of the inverse function theorem. Denoting the transformed function again by $f$, we obtain (\ref{eq:normal}).
Here, the coordinate transformation $\Psi$ is understood as the composition of all previous ones.
\qed

\begin{proposition}[Nondegenerate minimizers]
\label{prop:nondeg-min}
Let $\bar x$ be a nondegenerate M-stationary point for SCNO. Then, $\bar x$ is a local
minimizer for SCNO if and only if its M-index vanishes.
\end{proposition}

\proof
Let $\bar x$ be a nondegenerate M-stationary point for SCNO. The application of Morse Lemma from Theorem \ref{thm:morse} says that there exist neighborhoods $U_{\bar x}$ and $V_0$ of $\bar x$ and $0$, respectively, and a local $C^1$-coordinate system $\Psi: U_{\bar x} \rightarrow V_0$ of $\R^n$ around $\bar x$ such that:
\begin{equation}
\label{eq:help-nm0}
f\circ \Psi^{-1}(y)= f(\bar x) + \sum\limits_{i=1}^{n-k}y_i + \sum\limits_{j=n-k+1}^{n}\pm y_j^2,
\end{equation}
where $y \in \R^{n-k,s-k}\times \R^k$. Therefore, $\bar x$ is a local minimizer for SCNO if and only if $0$ is a local minimizer of $f\circ \Psi^{-1}$ on the set $\left(\R^{n-k,s-k}\times \R^k\right) \cap V_0$. 
If the M-index of $\bar x$ vanishes, we have $k=s$ and $QI=0$, and (\ref{eq:help-nm0}) reads as
\begin{equation}
\label{eq:help-nm1}
f\circ \Psi^{-1}(y)= f(\bar x) + \sum\limits_{j=n-s+1}^{n} y_j^2,    
\end{equation}
where $y \in \{0\}^{n-s}\times \R^s$. Thus, $0$ is a local minimizer for (\ref{eq:help-nm1}). Vice versa, if $0$ is a local minimizer for (\ref{eq:help-nm0}), then obviously $k=s$ and $QI=0$, hence, the M-index of $\bar x$ vanishes.
\qed


Let $C^2\left(\R^n,\R\right)$ be endowed with the strong (or Whitney) $C^2$-topology, denoted by $C^k_s$ (see e.\,g. \cite{hirsch:1976}). The $C^k_s$-topology is generated by allowing perturbations of the functions, their gradients and Hessians, which are controlled by means of continuous positive functions. We say that 
a set is $C^2_s$-generic if it contains a countable intersection of $C^2_s$-open and -dense subsets. Since $C^2\left(\R^n,\R\right)$ endowed with the $C^2_s$-topology is a Baire space, generic sets are in particular dense.

\begin{theorem}[Genericity for SCNO]
\label{thm:generic}
Let $\mathcal{F} \subset C^2(\R^{n},\R)$ denote the subset of 
objective functions in SCNO for which each 
M-stationary point is nondegenerate. Then, $\mathcal{F}$ is $C^2_s$-open and -dense.
\end{theorem}

\proof
Let us fix a number of non-zero entries $k \in \{0, \ldots,s\}$, an index set of $k$ non-zero entries $D \subset \{1,\ldots,n\}$, i.\,e. $|D|=k$, an index subset of zero entries $E \subset \{1,\ldots,n\} \backslash D$, and a rank $r \in \{0, \ldots, k\}$. For this choice we consider the set $\Gamma_{k,D,E,r}$ of $x$ such that the following conditions are satisfied:
\begin{itemize}
    \item [] (m1) $x_i\not =0$ for all $i \in D$, and  $x_i=0$ for all $i \in \{1,\ldots,n\} \backslash D$,
    \item [] (m2) $\displaystyle \frac{\partial  f}{\partial  x_i}(x) = 0$ for all $i \in D$,
    \item [] (m3) if $k < s$ then $\displaystyle \frac{\partial  f}{\partial  x_i}(x) = 0$ for all $i \in E$,
    \item [] (m4) the matrix $\displaystyle \left(\frac{\partial^2 f}{\partial x_i \partial x_j}\left(x \right) \right)_{i,j \in D}$ has rank $r$.
\end{itemize}
Note that (m1) refers to feasibility, (m2) to M-stationarity, and (m3)-(m4) describe possible violations of ND1-ND2, respectively.

Now, it suffices to show that all $\Gamma_{k,D,E,r}$ are generically empty whenever $E$ is nonempty or the rank $r$ is less than $k$. 
By setting $I_1(x)=D$ and $I_0(x)=\{1,\ldots,n\} \backslash D$, this would mean, respectively, that at least one of the derivatives $\displaystyle \frac{\partial f}{\partial x_i} \left( x \right)$  vanishes for $i\in E \subset I_0(x)$ in ND1 if $k <s$, or the matrix $\displaystyle \left(\frac{\partial^2 f}{\partial x_i \partial x_j}\left(x \right) \right)_{i,j \in I_1(x)}$ is singular in ND2.
In fact, the available degrees of freedom of the variables involved in each $\Gamma_{k,D,E,r}$ are $n$. The loss of freedom caused by (m1) is $n-k$, and the loss of freedom caused by (m2) is $k$. Hence, the total loss of freedom is $n$. We conclude that a further nondegeneracy would exceed the total available degrees of freedom $n$. By virtue of the jet transversality theorem from \cite{Jongen:2000}, generically the sets $\Gamma_{k,D,E,r}$ must be empty.

For the openness result, we argue in a standard way. Locally, M-stationarity can be written via stable equations. Then, the implicit function theorem for Banach spaces can be applied to follow M-stationary points with respect to (local) $C^2$-perturbations of defining functions. Finally, a standard globalization procedure exploiting the specific properties of the strong $C^2_s$-topology can be used to construct a (global) $C_s^2$-neighborhood of problem data for which the nondegeneracy property is stable.
\qed

\begin{theorem}[Genericity for minimizers]
Generically, all minimizers of SCNO are nondegenerate with the vanishing M-index. 
\end{theorem}

\proof Note that every local minimizer of SCNO has to be M-stationary. Nondegenerate M-stationary points are generic by Theorem \ref{thm:generic}. Hence, generically, local minimizers are nondegenerate. Due to Proposition \ref{prop:nondeg-min}, they have vanishing M-index. \qed

By recalling Definition \ref{def:mindex} of M-index, we deduce the following important Corollary \ref{cor:min} on the structure of minimizers for SCNO.

\begin{corollary}[Sparsity constraint at minimizers]
\label{cor:min}
At each generic local minimizer $\bar x \in \R^{n,s}$ of SCNO the sparsity constraint is active, i.\,e. $\left\|\bar x\right\|_0=s$.
\end{corollary}

\section{Relation to other stationarity concepts}
\label{sec:rel}

We relate M-stationarity to other well-known stationarity concepts for SCNO from the literature. First, we focus on S-stationarity introduced in \cite{burdakov:2016}. Then, the notions of basic feasibility and CW-minimality from \cite{beck:2013} will be discussed.

\subsection{S-stationarity}

In \cite{burdakov:2016} the following observation has been made: $\bar x$ solves SCNO if and only if there exists $\bar y$ such that $\left(\bar x, \bar y\right)$ solves the following mixed-integer program:
\begin{equation}
   \label{eq:mix-int}
   \min_{x,y} \,\, f(x) \quad \mbox{s.\,t.} \quad 
   \sum_{i=1}^{n} y_i \geq n - s,  \quad 
      y_i \in \{0,1\},  \quad 
      x_iy_i =0,  \quad  i=1, \ldots, n.
\end{equation}
Using the standard relaxation of the binary constraints $y_i\in\{0,1\}$, the authors arrive at the following continuous optimization problem:
\begin{equation}
   \label{eq:relax}
   \min_{x,y} \,\,f(x) \quad \mbox{s.\,t.} \quad 
   \sum_{i=1}^{n} y_i \geq n - s,  \quad 
      y_i \in [0,1],  \quad 
      x_iy_i =0,  \quad  i=1, \ldots, n.
\end{equation}
 As pointed out in \cite{burdakov:2016}, SCNO and the optimization problem (\ref{eq:relax}) are closely related: $\bar x$ solves SCNO if and only if there exists a vector $\bar y$ such that $\left(\bar x, \bar y\right)$ solves (\ref{eq:relax}). Additionally, the concept of  S-stationarity is proposed for (\ref{eq:relax}). For its formulation the following index sets are needed:
\[
 \begin{array}{lcl}
\displaystyle I_{\pm0}\left(\bar x, \bar y\right) &=& \displaystyle \left\{ i \in \{1,\ldots, n\} \,\left|\, \bar x_i \not = 0, \bar y_i =0 \right.\right\}, \\
\displaystyle  I_{00}\left(\bar x, \bar y\right) &=& \displaystyle  \left\{ i \in \{1,\ldots, n\} \,\left|\, \bar x_i = 0, \bar y_i =0 \right.\right\}.
 \end{array}{}
\]


\begin{definition}[S-stationarity, \cite{burdakov:2016}]
\label{def:m-s}
A feasible point $\left(\bar x, \bar y\right)$ of (\ref{eq:relax}) is called S-stationary 
if there exist real multipliers $\gamma_1, \ldots, \gamma_n$, such that
\begin{equation}
\label{eq:Sstat}
      \nabla f \left(\bar x\right) + \sum_{i}^{n} \gamma_i e_i =0, \quad \gamma_i =0 \mbox{ for all } i \in I_{\pm0}\left(\bar x, \bar y\right),
\end{equation}  
and, additionally, it holds:
\[
\gamma_i =0 \mbox{ for all } i \in I_{00}\left(\bar x, \bar y\right).
\]
\end{definition}


\begin{remark}[M-stationarity]
We point out that initially \cite{burdakov:2016} defined the concept of M-stationarity for the relaxed optimization problem (\ref{eq:relax}). Namely, a feasible point $\left(\bar x, \bar y\right)$ of (\ref{eq:relax}) is called M-stationary 
if just (\ref{eq:Sstat}) is valid. Due to the feasibility of $\left(\bar x, \bar y\right)$, we have $\bar y_i = 0$ if
$\bar x_i \not =0$ for all $i=1,\ldots, n$. Hence, it holds:
\[
I_{00}\left(\bar x, \bar y\right) = I_1\left(\bar x\right),
\]
and  M-stationarity is independent from the auxiliary variable $\bar y$. Thus, already in \cite{bucher:2018} it is sometimes said that a feasible point $\bar x$ of SCNO is M-stationary itself. We use M-stationarity exactly in this sense, cf. Definition \ref{def:m-stat}.
\end{remark}

In order to relate M- and S-stationarity, we introduce the canonical choice of the auxiliary variables $\bar y$ for a feasible point $\bar x$ of SCNO:
\begin{equation}
    \label{eq:canon}
   \bar y_i = \left\{\begin{array}{ll}
        0, & \mbox{if } i \in I_{1}\left(\bar x\right), \\
        1, &  \mbox{if } i \in I_{0}\left(\bar x\right).
   \end{array} \right.
\end{equation}
The auxiliary variables $\bar y$ can be seen as counters of the zero elements of $\bar x$. Note that $\left(\bar x, \bar y\right)$ becomes feasible for (\ref{eq:relax}).

\begin{proposition}[M- and S-stationarity]
\label{prop:t-m-s}
If $\left(\bar x, \bar y\right)$ is S-stationary for (\ref{eq:relax}) then $\bar x$ is M-stationary for SCNO. Vice versa, for any M-stationary point $\bar x$ the canonical choice (\ref{eq:canon}) of auxiliary variables $\bar y$ provides an S-stationary point $\left(\bar x, \bar y\right)$ for (\ref{eq:relax}). 
\end{proposition}

\proof Let $\left(\bar x, \bar y\right)$ be S-stationary for (\ref{eq:relax}).
After a moment of reflection we see that $I_{\pm0}\left(\bar x, \bar y\right) = I_{1}\left(\bar x\right)$ is the support of $\bar x$, and (\ref{eq:Sstat}) reads as the M-stationarity of $\bar x$:
\[
    \nabla_i f \left(\bar x\right) = 0 \mbox{ for all } i \in I_{1}\left(\bar x\right). 
\]

Vice versa, let $\bar x$ be an M-stationary point for SCNO with the canonical choice (\ref{eq:canon}) of $\bar y$. Then, $\left(\bar x, \bar y\right)$ is feasible for (\ref{eq:relax}), since 
\[
   \sum_{i=1}^{n} \bar y_i = \left| I_{0}\left(\bar x\right) \right| = n - \left| I_{1}\left(\bar x\right) \right| \geq n - s.
\]
The last inequality is due to $\left\|\bar x\right\|_0 \leq s$ or, equivalently, $\left| I_{1}\left(\bar x\right) \right| \leq s$. 
Moreover, by the choice of $\bar y$ we have $I_{\pm0}\left(\bar x, \bar y\right) = I_{1}\left(\bar x\right)$ and $I_{00}\left(\bar x, \bar y\right) = \emptyset$. Thus, due to the M-stationarity of $\bar x$, (\ref{eq:Sstat}) is fulfilled, and $\left(\bar x, \bar y\right)$ is S-stationary. \qed

The importance of $S$-stationary points is due to the following Proposition \ref{prop:s-kkt}.

\begin{proposition}[S-stationarity and KKT-points, \cite{burdakov:2016}]
\label{prop:s-kkt} 
A feasible point $\left(\bar x, \bar y\right)$ satisfies the Karush-Kuhn-Tucker condition if and only if it is S-stationary for (\ref{eq:relax}). 
\end{proposition}

Despite this appealing relation, nondegenerate M-stationary points of SCNO may cause degeneracies of the corresponding S-stationary points. This means that they become degenerate Karush-Kuhn-Tucker-points for (\ref{eq:relax}), i.\,e. the linear independent constraint qualification, strict complementarity, or second-order regularity is violated. The appearance of these degeneracies is mainly due to the fact that the objective function in (\ref{eq:relax}) does not depend on $y$-variables. We illustrate this phenomenon by means of the following Example \ref{ex:s-deg}.

\begin{example}[S-stationarity and degeneracies]
\label{ex:s-deg}
We consider the following SCNO with $n=2$ and $s=1$:
\[
   \min_{x_1,x_2}\,\, \left(x_1-1\right)^2 + \left(x_2-1\right)^2 \quad \mbox{s.\,t.} \quad 
   \left\|\left(x_1, x_2\right)\right\|_0 \leq 1.
\]
It is easy to see that the feasible point $\bar x=(0,0)$ is M-stationary with $\left\|\bar x\right\|_0=k=0$. Moreover, it is nondegenerate
with quadratic index $QI=0$. For its M-index we have 
\[
s-k+QI=1-0+0=1,
\]
meaning that $\bar x$ is a saddle point which connects two minimizers $(1,0)$ and $(0,1)$. Further, by the canonical choice (\ref{eq:canon}) of auxiliary $y$-variables, we obtain the corresponding S-stationary point $(\bar x, \bar y)=(0,0,1,1)$. Due to Proposition \ref{prop:s-kkt}, $(\bar x, \bar y)$ is also a Karush-Kuhn-Tucker-point for the optimization problem (\ref{eq:relax}):
\[
   \min_{x,y} \,\, \left(x_1-1\right)^2 + \left(x_2-1\right)^2 \quad \mbox{s.\,t.} \quad 
   y_1+y_2 \geq 1,  \quad 
      y_1, y_2 \in [0,1],  \quad 
      x_1y_1 =0,\quad  x_2y_2=0.
\]
The gradients of the active constraints at $(\bar x, \bar y)$ are linearly independent:
\[
  \left(\begin{array}{c}
        0 \\ 0 \\ 1 \\ 0
  \end{array}\right), 
  \left(\begin{array}{c}
        0 \\ 0 \\ 0 \\ 1
  \end{array}\right),
  \left(\begin{array}{c}
        \bar y_1 \\ 0 \\ \bar x_1 \\ 0
  \end{array}\right)=\left(\begin{array}{c}
        1 \\ 0 \\ 0 \\ 0
  \end{array}\right),
  \left(\begin{array}{c}
        0 \\ \bar y_2 \\ 0 \\ \bar x_2
  \end{array}\right)=\left(\begin{array}{c}
        0 \\ 1 \\ 0 \\ 0
  \end{array}\right).
\]
Hence, the linear independent constraint qualification holds at $(\bar x, \bar y)$.
Let us determine the unique Lagrange multipliers from the Karush-Kuhn-Tucker condition:
\[
   \left(\begin{array}{c}
        2 (\bar x_1 -1)\\ 2 (\bar x_2 -1) \\ 0 \\ 0
  \end{array}\right)= \mu_1 \left(\begin{array}{c}
        0 \\ 0 \\ 1 \\ 0
  \end{array}\right) + \mu_2 \left(\begin{array}{c}
        0 \\ 0 \\ 0 \\ 1
  \end{array}\right) + \lambda_1 \left(\begin{array}{c}
        1 \\ 0 \\ 0 \\ 0
  \end{array}\right) + \lambda_2 \left(\begin{array}{c}
        0 \\ 1 \\ 0 \\ 0
  \end{array}\right), \quad \mu_1, \mu_2 \leq 0.
\]
We get $\mu_1=\mu_2=0$ and $\lambda_1=\lambda_2=-2$. Hence, the strict complementarity is violated at $(\bar x, \bar y)$. Finally, the tangential space on the feasible set vanishes at $(\bar x, \bar y)$. Hence, the second derivative of the corresponding Lagrange function restricted to the tangential space is trivially nonsingular. This means that the second-order regularity is fulfilled at $(\bar x, \bar y)$. Overall, we claim that $(\bar x, \bar y)$ is a degenerate Karush-Kuhn-Tucker-point for  (\ref{eq:relax}) due to the lack of strict complementarity. It remains to note that the degeneracy of S-stationary points $\left(\bar x, y\right)$ prevails if other choices of auxiliary $y$-variables are made. \qed
\end{example}

An attempt to define a tailored notion of nondegeneracy for S-stationary points of (\ref{eq:relax}) has been recently undertaken in \cite{bucher:2018}. Let us briefly recall their main idea. For that, the so-called CC-linearization cone $\mathcal{L}^{CC}\left(\bar x, \bar y\right)$ at a feasible point $\left(\bar x, \bar y\right)$ of (\ref{eq:relax}) is used, cf. \cite{cervinka:2016}. Namely,
\[
   \left(d_x,d_y\right) \in \mathcal{L}^{CC}\left(\bar x, \bar y\right) \subset \R^n \times \R^n  
\]
satisfies by definition the following conditions:
\begin{equation}
\label{eq:cc-l}
\left\{
\begin{array}{l}
    \displaystyle \sum_{i=1}^{n} \left(d_y\right)_i \geq 0 \mbox{ if } \displaystyle \sum_{i=1}^{n} \bar y_i = n-s, \\ 
     \left(d_y\right)_i=0 \mbox{ for all } i \in I_{\pm0}\left(\bar x, \bar y\right), \\ 
     \left(d_y\right)_i\geq 0 \mbox{ for all } i \in I_{00}\left(\bar x, \bar y\right), \\ 
     \left(d_y\right)_i\leq 0 \mbox{ for all } i \in I_{01}\left(\bar x, \bar y\right), \\
    \left(d_x\right)_i= 0 \mbox{ for all } i \in I_{01}\left(\bar x, \bar y\right) \cup I_{0+}\left(\bar x, \bar y\right), \\
    \left(d_x\right)_i \left(d_y\right)_i= 0 \mbox{ for all } i \in I_{00}\left(\bar x, \bar y\right).
\end{array}
\right.   
\end{equation}
Here, the new index sets are
\[
 \begin{array}{lcl}
\displaystyle I_{01}\left(\bar x, \bar y\right) &=& \displaystyle \left\{ i \in \{1,\ldots, n\} \,\left|\, \bar x_i = 0, \bar y_i =1 \right.\right\}, \\
\displaystyle  I_{0+}\left(\bar x, \bar y\right) &=& \displaystyle  \left\{ i \in \{1,\ldots, n\} \,\left|\, \bar x_i = 0, \bar y_i \in (0,1) \right.\right\}.
 \end{array}{}
\]

\begin{definition}[CC-SOSC, \cite{bucher:2018}]
\label{def:s-sosc}
Let $\left(\bar x, \bar y\right)$ be an S-stationary point for (\ref{eq:relax}). If for all directions $\left(d_x,d_y\right) \in \mathcal{L}^{CC}\left(\bar x, \bar y\right)$ with $d_x \not =0$, we have
\[
  d_x^T \cdot D^2f(\bar x) \cdot d_x > 0,
\]
then the cardinality constrained second-order sufficient optimality condition (CC-SOSC) is said to hold at $\left(\bar x, \bar y\right)$.
\end{definition}
 
The role of CC-SOSC can be seen from the following Proposition \ref{prop:min-sosc}.

\begin{proposition}[Sufficient optimality condition, \cite{bucher:2018}]
\label{prop:min-sosc}
Let $\left(\bar x, \bar y\right)$ be an S-stationary point for (\ref{eq:relax}) satisfying CC-SOSC. Then, $\left(\bar x, \bar y\right)$ is a strict local minimizer of (\ref{eq:relax}) with respect to $x$, i.\,e.
\[
  f\left(\bar x\right) < f(x)
\]
for all feasible points $(x,y)$ of (\ref{eq:relax}) taken sufficiently close to $\left(\bar x, \bar y\right)$, and fulfilling $x \not = \bar x$. 
\end{proposition}

We relate the concepts of nondegeneracy for M-stationary points and of CC-SOSC for S-stationary points.

\begin{proposition}[Nondegeneracy and CC-SOSC]
\label{prop:nondeg-sosc}
Let $\bar x$ be an M-stationary point for SCNO with $\left\|\bar x\right\|_0=s$.
Assume that CC-SOSC holds at the S-stationary point $\left(\bar x, \bar y\right)$ for (\ref{eq:relax}) with the canonical choice (\ref{eq:canon}) of auxiliary variables $\bar y$. Then, $\bar x$ is a nondegenerate local minimizer for SCNO. 
\end{proposition}

\proof Due to the canonical choice (\ref{eq:canon}) of auxiliary variables $\bar y$, the index sets from the definition of the CC-linearization cone $\mathcal{L}^{CC}\left(\bar x, \bar y\right)$ are
\[
     I_{\pm0}\left(\bar x, \bar y\right)= I_1\left(\bar x\right), \quad
     I_{00}\left(\bar x, \bar y\right) = I_{0+}\left(\bar x, \bar y\right) =\emptyset, \quad
     I_{01}\left(\bar x, \bar y\right)=I_0\left( \bar x\right).
\]
Due to $\left\|\bar x\right\|_0=s$, we additionally have $\displaystyle \sum_{i=1}^{n} \bar y_i = n-s$. Recalling (\ref{eq:cc-l}), $\left(d_x,d_y\right) \in \mathcal{L}^{CC}\left(\bar x, \bar y\right)$ if and only if
\[
\left\{
\begin{array}{l}
  \displaystyle \sum_{i=1}^{n} \left(d_y\right)_i \geq 0, \\
     \left(d_y\right)_i=0 \mbox{ for all } i \in I_{1}\left(\bar x\right), \\ 
     \left(d_y\right)_i\leq 0 \mbox{ for all } i \in I_{0}\left(\bar x\right), \\
    \left(d_x\right)_i= 0 \mbox{ for all } i \in I_{0}\left(\bar x\right).   
\end{array}
\right.
\]
Hence, it holds:
\[
 \mathcal{L}^{CC}\left(\bar x, \bar y\right) = \left\{
 \left(d_x,0\right) \, \left|\,
    \left(d_x\right)_i= 0 \mbox{ for all } i \in I_{0}\left(\bar x\right)
 \right.
 \right\},
\]
so that CC-SOSC says that the matrix $\displaystyle \left(\frac{\partial^2 f}{\partial x_i \partial x_j}\left(\bar x \right) \right)_{i,j \in I_1\left(\bar x\right)}$ is positive definite. By Definition \ref{def:nondeg}, the M-stationary point $\bar x$ is then nondegenerate with the vanishing quadratic index, i.\,e. $QI=0$. Thus, using again that $\left\|\bar x\right\|_0=s$, its M-index becomes $s-s+QI=0$. Finally, Proposition \ref{prop:nondeg-min} provides the assertion. \qed 

If the sparsity constraint is not active for an M-stationary point $\bar x$ of SCNO, i.\,e. $\left\|\bar x\right\|_0<s$, the implication in Proposition \ref{prop:nondeg-sosc} does not hold in general anymore. Namely, $\bar x$ does not need to be a local minimizer for SCNO, even if 
CC-SOSC holds at the corresponding S-stationary point $\left(\bar x, \bar y\right)$ with the canonical choice (\ref{eq:canon}) of auxiliary variables $\bar y$. This is illustrated by means of the following Example \ref{ex:sosc-binding}.

\begin{example}[Sparsity constraint and CC-SOSC]
\label{ex:sosc-binding}
We consider the following SCNO with $n=2$ and $s=1$:
\[
   \min_{x_1,x_2}\,\, x_1 + x_2 \quad \mbox{s.\,t.} \quad 
   \left\|\left(x_1, x_2\right)\right\|_0 \leq 1.
\]
It is easy to see that the feasible point $\bar x=(0,0)$ is M-stationary. 
Note that the sparsity constraint is not active for $\bar x$, since
$k=\left\|\bar x\right\|_0=0 < 1 =s$. By the canonical choice (\ref{eq:canon}) of auxiliary $y$-variables, we obtain the corresponding S-stationary point $(\bar x, \bar y)=(0,0,1,1)$. Analogously to the proof of Proposition \ref{prop:min-sosc} and by recalling (\ref{eq:cc-l}), $\left(d_x,d_y\right) \in \mathcal{L}^{CC}\left(\bar x, \bar y\right)$ if and only if
\[
\left\{
\begin{array}{l}
     \left(d_y\right)_i=0 \mbox{ for all } i \in I_{1}\left(\bar x\right), \\ 
     \left(d_y\right)_i\leq 0 \mbox{ for all } i \in I_{0}\left(\bar x\right), \\
     \left(d_x\right)_i= 0 \mbox{ for all } i \in I_{0}\left(\bar x\right).   
\end{array}
\right.
\]
Note that here $\displaystyle I_{1}\left(\bar x\right)=\emptyset$ and $\displaystyle I_{0}\left(\bar x\right)=\{1,2\}$. Hence, the CC-linearization cone is
\[
 \mathcal{L}^{CC}\left(\bar x, \bar y\right) = \left\{
 \left(0,d_y\right) \, \left|\,
    \left(d_y\right)_1, \left(d_y\right)_2\leq 0
 \right.
 \right\}.
\]
Overall, CC-SOSC trivially holds at $\left(\bar x, \bar y\right)$, and as follows from Proposition \ref{prop:min-sosc}, it is a strict local minimizer of (\ref{eq:relax}) with respect to $x$. Nevertheless, $\bar x$ is not a local minimizer. Actually, it is a nondegenerate M-stationary point with the quadratic index $QI=0$. For its M-index we have 
\[
s-k+QI=1-0+0=1.
\]
We conclude that $\bar x$ is rather a saddle point for SCNO. \qed
\end{example}

\subsection{Basic feasibility and CW-minimality}

We proceed by discussing stationarity concepts from \cite{beck:2013}.
Inspired by linear programming terminology, they first introduce the notion of a basic feasible vector for SCNO.

\begin{definition}[Basic feasibility, \cite{beck:2013}]
\label{def:bf}
A vector $\bar x \in \R^{n,s}$ with $\left\|\bar x \right\|_0=k$ is called basic feasible (BF) for SCNO if the following conditions are fulfilled: 
\begin{itemize}
    \item[] BF1: in case $k < s$, it holds: 
\[
\frac{\partial f}{\partial  x_i} \left( \bar x\right) = 0 \mbox{ for all } i =1,\ldots,n,
\]
    \item[] BF2: in case $k = s$, it holds: 
\[
\frac{\partial f}{\partial  x_i} \left( \bar x\right) = 0 \mbox{ for all } i \in I_1\left(\bar x\right).
\]
\end{itemize}
\end{definition}

Attention has been also paid to the notion of coordinate-wise minimum for SCNO.   

\begin{definition}[CW-minimality, \cite{beck:2013}]
\label{def:cw}
A vector $\bar x \in \R^{n,s}$ with $\left\|\bar x \right\|_0=k$ is called coordinate-wise (CW) minimum for SCNO if the following conditions are fulfilled: 
\begin{itemize}
    \item[] CW1: in case $k < s$, it holds: 
\[
f\left(\bar x \right)= \min_{t \in \R} f\left( \bar x +t e_i\right) \mbox{ for all } i =1,\ldots,n,
\]
    \item[] CW2: in case $k = s$, it holds: 
\[
f\left(\bar x \right) \leq \min_{t \in \R} f\left( \bar x - \bar x_i e_i+t e_j\right)  \mbox{ for all } i \in I_1\left(\bar x\right) \mbox { and } j=1,\ldots, n.
\]
\end{itemize}
\end{definition}

Basic feasibility and CW-minimality can be viewed as necessary optimality condition for SCNO.

\begin{proposition}[BF-vector and CW-minimum, \cite{beck:2013}]
\label{prop:bf-cw}
Every global minimizer for SCNO is a CW-minimum, and every CW-minimum for SCNO is a BF-vector.
\end{proposition}

It is claimed in \cite{beck:2013} that the basic feasibility condition is quite weak, namely, there are many BF-points that are not optimal for SCNO.
The notion of CW-minimum provides a much stricter necessary optimality condition. Based on the latter, a greedy sparse-simplex method for the numerical treatment of SCNO is proposed by \cite{beck:2013}. 
Let us now examine the relation between M-stationarity, basic feasibility, and  CW-minimality.

\begin{proposition}[M-stationarity, BF-vector, and CW-minimum]
\label{prop:t-bf-cw}
Every BF-vector for SCNO is an M-stationary point, in particular, so is every CW-minimum.
\end{proposition}

\proof Let $\bar x$ be a BF-vector for SCNO with $\left\|\bar x \right\|_0=k$.
If $k < s$, then BF1 implies M-stationarity of $\bar x$. If $k=s$, then BF2 coincides with the latter property. Since every CW-minimum for SCNO is a BF-vector according to Proposition \ref{prop:bf-cw}, the assertion follows. \qed

Proposition \ref{prop:t-bf-cw} says that M-stationarity is an even weaker condition than basic feasibility and CW-minimality. Why should we care about M-stationarity then? Is it not enough to rather focus on the stricter necessary optimality condition of CW-minimality as in \cite{beck:2013}? It turns out that CW-minima need not to be stable with respect to data perturbations. Namely, after an arbitrarily small $C^2$-perturbation of $f$ a locally unique CW-minimum may bifurcate into multiple CW-minima. More importantly, this bifurcation unavoidably causes the emergence of M-stationary points, being different from CW-minima. Next Example \ref{ex:cw-inst} illustrates this instability phenomenon.

\begin{example}[CW-mimimum and instability]
\label{ex:cw-inst}
We consider the following SCNO with $n=2$ and $s=1$:
\begin{equation}
\label{eq:inst}
\min_{x_1,x_2}\,\, x_1^2 + x_2^2 \quad \mbox{s.\,t.} \quad 
   \left\|\left(x_1, x_2\right)\right\|_0 \leq 1.
\end{equation}
Obviously, $\bar x=(0,0)$ is the unique minimizer of (\ref{eq:inst}). Due to Proposition \ref{prop:bf-cw}, it is also a CW-minimum, as well as a BF-vector. Further, let us perturb (\ref{eq:inst}) by using an arbitrarily small $\varepsilon >0$ as follows:  
\begin{equation}
\label{eq:st}
\min_{x_1,x_2}\,\, \left(x_1-\varepsilon\right)^2 + \left(x_2-\varepsilon\right)^2 \quad \mbox{s.\,t.} \quad 
   \left\|\left(x_1, x_2\right)\right\|_0 \leq 1.
\end{equation}
It is easy to see that the perturbed problem (\ref{eq:st}) has now two solutions $\bar x^1=(\varepsilon,0)$ and $\bar x^2=(0,\varepsilon)$. 
Both are CW-minima, and, hence, BF-points.
Here, we observe a bifurcation of the CW-minimum $\bar x$ of the original problem (\ref{eq:inst}) into two CW-minima $\bar x^1$ and $\bar x^2$ of the perturbed problem (\ref{eq:st}). Let us explain this bifurcation in terms of M-stationarity. The bifurcation is caused by the degeneracy of $\bar x$ viewed as an M-stationary point of the original problem (\ref{eq:inst}). Note that ND1 is violated at the M-stationary point $\bar x$ of the original problem (\ref{eq:inst}). More interestingly, 
although $\bar x$ is neither a CW-minimum nor a BF-vector of (\ref{eq:st}) anymore, it becomes a new M-stationary point for the perturbed problem. In fact, due to $\left\|\bar x\right\|_0=k=0$ and the validity of ND1, $\bar x$ is a nondegenerate M-stationary point of (\ref{eq:st})  
with the quadratic index $QI=0$. For its M-index we have 
\[
s-k+QI=1-0+0=1,
\]
meaning that $\bar x$ is a saddle point which connects two nondegenerate minimizers $\bar x^1$ and $\bar x^2$ of (\ref{eq:st}). Overall, we conclude that the degenerate CW-minimum $\bar x$ of the original problem (\ref{eq:inst}) is not stable. Moreover, it bifurcates into two nondegenerate CW-minima $\bar x^1$ and $\bar x^2$, as well as leads to one nondegenerate saddle point $\bar x$ of the perturbed problem (\ref{eq:inst}). 
\qed
\end{example}

Example \ref{ex:cw-inst} suggests to consider nondegenerate BF-vectors or nondegenerate CW-minima for SCNO, in order to guarantee their stability  with respect to sufficiently small data perturbations. Then, however, the sparsity constraint turns out to be active. This means that BF1 in Definition \ref{def:bf} and CW1 in Definition \ref{def:cw} become redundant. 

\begin{proposition}[BF-vector, CW-minumum and nondegeneracy]
\label{prop:bf-nondeg}
Let $\bar x$ be a BF-vector for SCNO with $\left\|\bar x\right\|_0=k$. If it is nondegenerate as an M-stationary point for SCNO, then $k=s$. The same applies for CW-minima.
\end{proposition}

\proof Assume that $k < s$, then ND1 contradicts BF1, whenever $I_0\left(\bar x\right) \not = \emptyset$. Otherwise, we have $k=n$, and, hence, $n < s$, a contradiction. It remains to note that every CW-minimum for SCNO is a BF-vector due to Proposition \ref{prop:bf-cw}. \qed

\section{Global results}
\label{sec:crit}

Let us study the topological properties of lower level sets
\[
M^{a}=\left\{ x \in\R^{n,s}\, \left\vert \, f(x)\le a \right. \right\},
\]
where $a \in \R$ is varying. For that, we define intermediate sets for $a<b$:
\[
M^{b}_{a}=\left\{ x \in\R^{n,s}\, \left\vert \, a \leq f(x) \leq b \right. \right\}.
\]
For the topological concepts used below we
refer to \cite{spanier:1966}.

Let us start with Assumption \ref{ass:prop} which is usual within the scope of Morse theory, cf. \cite{goresky:1988}. It prevents from considering asymptotic effects at infinity.

\begin{assumption}
\label{ass:prop}
The restriction of the objective function $f_{|\R^{n,s}}$  on the SCNO feasible set is proper, i.\,e. $f^{-1}(K)\cap \R^{n,s}$ is compact for any compact set $K \subset \R$.
\end{assumption}

\begin{theorem}[Deformation for SCNO]
\label{thm:def}
Let Assumption \ref{ass:prop} be fulfilled and $M^b_a$ contain no M-stationary points for SCNO. Then, $M^a$ is homeomorphic to $M^b$.
\end{theorem}

\proof
We apply Proposition 3.2 from Part I in \cite{goresky:1988}. The latter provides the deformation for general Whitney stratified sets with respect to critical points of proper maps. Note that the SCNO feasible set admits a Whitney stratification:
\[
   \R^{n,s} = \bigcup_{\scriptsize \begin{array}{c}
        I \subset \{1, \ldots,n\} \\
        |I| \leq s 
   \end{array}} \bigcup_{J \subset I} Z_{I,J},
\]   
where 
\[
  Z_{I,J} = \left\{x\in \R^n \,\left|\, x_{I^c} =0, x_J > 0, x_{I\backslash J} <0 \right.\right\}.
\]
The notion of criticality used in \cite{goresky:1988} can be stated for SCNO as follows. A point $\bar x \in \R^{n,s}$ is called critical for $f_{|\R^{n,s}}$ if it holds:
\[
\nabla f\left(\bar x\right)_{|T_{\bar x} Z} =0,
\]
where $Z$ is the stratum of $\R^{n,s}$ which contains $\bar x$, and $T_{\bar x} Z$ is the tangent space of $Z$ at $\bar x$. By identifying $I =I_1\left(\bar x\right)$ and, hence, $I^c =I_0\left(\bar x\right)$, we see that the concepts of criticality and M-stationarity coincide.
This concludes the assertion. \qed

Let us now turn our attention to the topological changes of lower level sets when passing an M-stationary level. Traditionally, they are described by means of the so-called cell-attachment.
We first consider a special case of cell-attachment. For that, let $N^\epsilon$ denote the lower level set of a special linear function on $\R^{p,q}$, i.\,e.
\[
 N^\epsilon = \left\{ x \in \R^{p,q} \,\left|\, \sum_{i=1}^{p} x_i \leq \epsilon \right. \right\},
\]
where $\epsilon \in \R$, and the integers $q < p$ are nonnegative. 

\begin{lemma}[Normal Morse data]
\label{lem:cat}
For any $\epsilon > 0$ the set $N^\epsilon$ is homotopy-equivalent to $N^{-\epsilon}$ with $\binom{p-1}{q}$ cells of dimension $q$ attached. The latter cells are the $q$-dimensional simplices from the collection 
\[
 \left\{ \left. \mbox{conv} \left(e_j, j \in J\right)  \,\right|\, J \subset \{1, \ldots,p\}, 1 \in J, |J| = q+1 \right\}.
\]
\end{lemma}

\proof Let $N_\epsilon$ denote the upper level set of a special linear function on $\R^{p,q}$, i.\,e.
\[
 N_\epsilon = \left\{ x \in \R^{p,q} \,\left|\, \sum_{i=1}^{p} x_i \geq \epsilon \right. \right\}.
\]
In terms of upper level sets Lemma \ref{lem:cat} can be obviously reformulated as follows: For any $\epsilon > 0$ the set $N_{-\epsilon}$ is homotopy-equivalent to $N_{\epsilon}$ with $\binom{p-1}{q}$ cells of dimension $q$ attached. Let us show the latter assertion.

First, we note that the sets $N_0$ and $N_{-\epsilon}$ are contractible. The contraction is performed via the mapping 
\[
   (x,t) \mapsto (1-t)\cdot x, \quad t \in [0,1].
\]
For the lower level set $N^{\epsilon}$ we have the representation
\[
   N^{\epsilon} = \bigcup\limits_{
\scriptsize \begin{array}{c}
        J \subset \{1, \ldots,p\} \\
        |J| = q 
   \end{array}} N^{\epsilon,J},
\]
where
\[
  N^{\epsilon,J} = \left\{ x \in \R^{p,q} \,\left|\, x_{J^c} =0, \sum_{i\in J} x_i \geq \epsilon \right. \right\}.
\]
Note that $N^{\epsilon,J}$ is homotopy-equivalent to the set $N^J$, where
\[
   N^J = \left\{ x \in \R^{p,q} \,\left|\, x_{J^c} =0, \sum_{i\in J} x_i = 1 \right. \right\}
\]
is the $(|J|-1)$-dimensional simplex $\mbox{conv} \left(e_j, j \in J\right)$ of $\R^p$. In fact, the map  
\[
(x,t) \mapsto t \cdot \frac{x}{\displaystyle \sum_{i= 1}^{p} x_i}+(1-t) \cdot x, \quad t \in [0,1]
\]
can be used for all $N^J$.
Altogether, $N_{\epsilon}$ is homotopy-equivalent to
\begin{equation}
    \label{eq:skel}
    \bigcup\limits_{
\scriptsize \begin{array}{c}
        J \subset \{1, \ldots,p\} \\
        |J| = q 
   \end{array}} \mbox{conv}\left(e_j, j \in J\right).
\end{equation}
Note that the set in (\ref{eq:skel}) is the $(q-1)$-skeleton of the $(p-1)$-dimensional simplex of $\R^p$. The $(q-1)$-skeleton of the $(p-1)$-dimensional simplex is the union of its simplices up to dimension $q-1$, see e.\,g. \cite{goerss:2009}.

Within the $(q-1)$-skeleton (\ref{eq:skel}), we close all $q$-dimensional holes by  attaching $q$-dimensional cells from the collection of simplices
\[
 \left\{ \left. \mbox{conv} \left(e_j, j \in J\right)  \,\right|\, J \subset \{1, \ldots,p\}, |J| = q+1 \right\}.
\]
The attachment should result in a contractible set, as it is actually $N_0$.
We note that the union of the subdivision 
\begin{equation}
    \label{eq:subdiv}
 \left\{ \left. \mbox{conv} \left(e_j, j \in J\right)  \,\right|\, J \subset \{1, \ldots,p\}, 1 \in J, |J| = q+1 \right\}
\end{equation}
is also contractible, namely, to $e_1$. To see this, we may use the map 
\[
(x,t) \mapsto t \cdot e_1+(1-t) \cdot x, \quad t \in [0,1].
\]
Furthermore, none of the relative interiors of the simplices in (\ref{eq:subdiv}) can be deleted. In fact, deleting gives rise to the boundary of a $q$-dimensional simplex and the latter is not contractible.

On the other hand, for any $J^* \subset \{1, \ldots,p\} \backslash \{1\}$ with $\left|J^*\right|=q+1$ the union
\begin{equation}
    \label{eq:gegen}
  \mbox{conv} \left(e_j, j \in J^*\right) \cup 
   \bigcup\limits_{
\scriptsize \begin{array}{c}
        J^{**} \subset J^* \\
        \left|J^{**}\right| = q 
   \end{array}} \mbox{conv}\left(e_j, j \in J^{**}\cup \{1\}\right)
\end{equation}
forms the boundary of the $(q+1)$-dimensional simplex $\mbox{conv} \left(e_j, j \in J^* \cup \{1\}\right)$. Hence, the set in (\ref{eq:gegen}) is not contractible.
Altogether, precisely the $q$-dimensional cells in (\ref{eq:subdiv}) can be attached to the $(q-1)$-skeleton (\ref{eq:skel}) in order to obtain a contractible set. Its number obviously equals 
$\binom{p-1}{q}$. This completes the proof. \qed

\begin{theorem}[Cell-Attachment for SCNO]
\label{thm:cell-a}
Let Assumption \ref{ass:prop} be fulfilled and
$M^b_a$ contain exactly one M-stationary point $\bar x$ for SCNO with $\left\|\bar x\right\|_0=k$ and the M-index equal to $s-k+QI$. If $a<f\left(\bar x \right) <b$,
then $M^b$ is homotopy-equivalent to $M^a$ with $\binom{n-k-1}{s-k}$ cells of dimension $s-k+QI$ attached, namely:
\[
 \bigcup_{\scriptsize \begin{array}{c}
        J \subset \{1, \ldots,n-k\} \\
        1 \in J, |J| = s-k +1
   \end{array}} \mbox{conv} \left(e_j, j \in J\right) \times [0,1]^{QI}.
\]
\end{theorem}

\proof
Theorem \ref{thm:def} allows deformations up to an arbitrarily small neighborhood of the M-stationary point $\bar x$. In such a neighborhood, we may assume without loss of generality that $\bar x=0$ and $f$ has the following form as from Theorem \ref{thm:morse}:
\begin{equation}
    \label{eq:n1}
        f(x)= f\left(\bar x\right) + \sum\limits_{i=1}^{n-k}x_i + \sum\limits_{j=n-k+1}^{n} \pm x_j^2,
\end{equation}
where $x \in \R^{n-k,s-k}\times \R^k$, and the number of negative squares in (\ref{eq:n1}) equals $QI$.
In terms of \cite{goresky:1988} the set 
$\R^{n-k,s-k}\times \R^k$ can be interpreted as the product of the tangential part $\R^k$ and the normal part $\R^{n-k,s-k}$. The cell-attachment along the tangential part is standard.  Analogously to the unconstrained case, one $QI$-dimensional cell has to be attached on $\R^k$.  The cell-attachment along the normal part is more involved. Due to Lemma \ref{lem:cat}, we need to attach $\binom{n-k-1}{s-k}$ cells on $\R^{n-k,s-k}$, each of dimension $s-k$. Finally, we apply Theorem 3.7 from Part I in \cite{goresky:1988}, which says that the local Morse data is the product of tangential and normal Morse data. Hence, the dimensions of the attached cells add together. Here, we have then to attach $\binom{n-k-1}{s-k}$ cells on $\R^{n-k,s-k}\times \R^k$, each of dimension $s-k+QI$. \qed

Let us put Theorem \ref{thm:cell-a} into the context of Morse theory as developed in the literature for other nonsmooth optimization problems. The new issue for SCNO is the multiplicity of attached cells. 

\begin{remark}[Multiplicity of attached cells]
We recall that for nonlinear programming problems (NLP)
the dimension of the cell to be attached while passing a critical point equals to its quadratic index, see e.\,g. \cite{Jongen:2000}. The situation changes if we consider mathematical programs with complementarity constraints (MPCC). Here, the dimension of attached cells equals to the  so-called C-index of C-stationary points, see \cite{shikhman:2009}. In addition to quadratic, the C-index also has a bi-active part. The latter counts negative pairs of Lagrange multipliers corresponding to the bi-active complementarity constraints. The cell-attachment for mathematical programs with vanishing constraints (MPVC) is similar, see \cite{shikhman:2012a}. The dimension of attached cells equals here to the so-called T-index of T-stationary points. The T-index consists again of quadratic and bi-active parts. We emphasize that the cell-attachment for SCNO considerably differs from the described cases of NLP, MPCC, and MPVC. The main difference is that {\it multiple} cells are involved into the cell-attachment procedure for SCNO. The multiplicity of attached cells is a novel and striking phenomenon in nonsmooth optimization not observed in the literature before. From the technical point of view, this makes the cell-attachment result for SCNO to appear rather challenging. Note that the determination of the number of attached cells becomes an involved combinatorial problem from algebraic topology, see Lemma \ref{lem:cat}.
\end{remark}

Let us present a global interpretation of our results for SCNO. For that, we need to state another assumption. Following Assumption \ref{ass:lb} is standard in the context of SCNO, cf. \cite{beck:2013}, and gives a necessary condition for its solvability.

\begin{assumption}
\label{ass:lb}
The restriction of the objective function $f_{|\R^{n,s}}$ on the SCNO feasible set is lower bounded.
\end{assumption}

Now, we consider M-stationary points $\bar x$ for SCNO with $\left\|\bar x\right\|_0=k$ and the M-index equal to one, thus, fulfilling $s-k+QI=1$. These so-called saddle points can be of two types:
\begin{itemize}
    \item[(I)] with active sparsity constraint and quadratic index equal to one, i.\,e. 
    \[
       k=s, \quad QI=1,
    \]
    \item[(II)] with exactly $s-1$ non-zero entries and vanishing quadratic index, i.\,e. 
    \[
       k=s-1, \quad QI=0.
    \]   
\end{itemize}

\begin{theorem}[Morse relation for SCNO]
\label{thm:mrel}
Let Assumptions \ref{ass:prop} and \ref{ass:lb} be fulfilled, and all M-stationary points of SCNO be nondegenerate. 
Additionally, we assume that there exists a connected lower level set which contains all M-stationary points. Then, it holds:
\begin{equation}
    \label{eq:mr}
    r_I + (n-s)r_{II} \geq r-1,
\end{equation}
where $r$ is the number of local minimizers of SCNO, $r_I$ and $r_{II}$ are the numbers of M-stationary points with M-index equal to one, which correspond to the types (I) and (II), respectively.
\end{theorem}

\proof 
We assume without loss of generality that the objective function $f$ has pairwise different values at all M-stationarity points of SCNO. If it is not the case, we may enforce this property by sufficiently small perturbations of the objective function.  
Due to the openness part in Theorem \ref{thm:generic}, all M-stationarity points of such a perturbed SCNO remain nondegenerate. Moreover, the formula (\ref{eq:mr}) is still valid since it does not depend on the functional values of $f$.

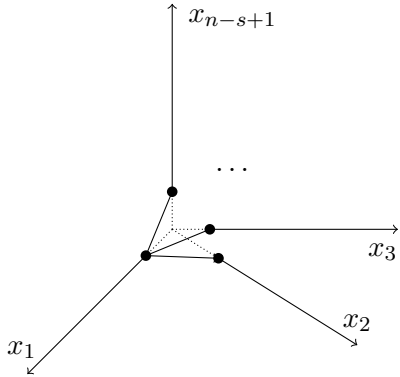
\begin{figure}
    \centering
\begin{tikzpicture}
\draw[->] (0.5,0) -- (xyz cs:x=3);
\node at (2.8,-0.3){$x_3$};
\draw[->] (0,0.5) -- (xyz cs:y=3);
\node at (0.8,2.8){$x_{n-s+1}$};
\draw[->] (-0.35,-0.35) -- (xyz cs:z=5);
\node at (-2.0,-1.6){$x_1$};
\node at (0.8,0.8){$\ldots$};
\draw[] (-0.35,-0.35) -- (0.5,0);
\draw[] (-0.35,-0.35) -- (0,0.5);
\draw[densely dotted](0,0) -- (xyz cs:y=0.5);
\draw[densely dotted](0,0) -- (xyz cs:x=0.5);
\draw[densely dotted](0,0) -- (xyz cs:z=1);
\draw[densely dotted](0,0) -- (xyz cs:x=1, z=1);
\draw[->] (1,0,1) -- (4,0,4);
\node at (4,0.3,4){$x_2$};
\draw[->] (-0.35,-0.35) -- (1,0,1);
\fill (0.5,0) circle(2pt);
\fill (0,0.5) circle(2pt);
\fill (-0.35,-0.35) circle(2pt);
\fill (1,0,1) circle(2pt);
   \end{tikzpicture} 
   \caption{Cell-attachment for type (II)}
    \label{fig:2}
\end{figure}

Further, let $q_a$ denote the number of connected components of the lower level set $M^a$.
We focus on how $q_a$ changes as $a \in \R$ increases. Due to Theorem \ref{thm:def}, $q_a$ can change only if passing through a value corresponding to an M-stationary point $\bar x$, i.\,e. $a=f\left(\bar x\right)$. 
In fact, Theorem \ref{thm:def} allows homeomorphic deformations of lower level sets up to an arbitrarily small neighborhood of the M-stationary point $\bar x$.
Then, we have to estimate the difference between $q_a$ and $q_{a-\varepsilon}$, where $\varepsilon > 0$ is arbitrarily, but sufficiently small, and $a=f\left(\bar x\right)$. This is done by a {\it local argument}. For that, let the M-index of $\bar x$ be $s-k+QI$ with $\left\|\bar x\right\|_0=k$. We use Theorem \ref{thm:cell-a} which says that $M^{a}$ is homotopy-equivalent to $M^{a-\varepsilon}$ with a cell-attachment of
\begin{equation}
    \label{eq:hep4}
 \bigcup_{\scriptsize \begin{array}{c}
        J \subset \{1, \ldots,n-k\} \\
        1 \in J, |J| = s-k +1
   \end{array}} \mbox{conv} \left(e_j, j \in J\right) \times [0,1]^{QI}.
\end{equation}
Let us distinguish the following cases:
\begin{itemize}
    \item[1)] $\bar x$ is a local minimizer with vanishing M-index, i.\,e. $k=s$ and $QI=0$.
    Then, by (\ref{eq:hep4}) we attach to $M^{a-\varepsilon}$ the cell $\mbox{conv}\left(e_1\right)$ of dimension zero. Consequently, a new connected component is created, and it holds: 
\[
   q_a = q_{a-\varepsilon} + 1.
\]
    \item[2)] $\bar x$ is of type (I) with M-index equal to one, i.\,e. $k=s$ and $QI=1$.
    Then, by (\ref{eq:hep4}) we attach to $M^{a-\varepsilon}$ the cell $\mbox{conv}\left(e_1\right)\times [0,1]$ of dimension one.
    Consequently, at most one connected component disappears, and it holds:
\[
   q_{a-\varepsilon} -1 \leq q_a \leq q_{a-\varepsilon}.
\]
 This case is well known from nonlinear programming, see e.\,g. \cite{Jongen:2000}.
    \item[3)] $\bar x$ is of type (II) with M-index equal to one, i.\,e. $k=s-1$ and $QI=0$.
    Then, by (\ref{eq:hep4}) we attach to $M^{a-\varepsilon}$ as many as $n-s+1$ cells of dimension one, namely:
\[
   \bigcup_{\scriptsize \begin{array}{c}
        j=2, \ldots, n-s+1 
   \end{array}} \mbox{conv} \left(e_1,e_j\right).
\]
     Consequently, at most $n-s$ connected components disappear, and it holds:
\[
   q_{a-\varepsilon} - (n-s) \leq q_a \leq q_{a-\varepsilon}.
\]
For illustration we refer to Figure \ref{fig:2}. Case 3) is new and characteristic for SCNO.
    \item[4)] $\bar x$ is M-stationary with M-index greater than one, i.\,e. $s-k+QI > 1$.
    The boundary of the cell-attachment in (\ref{eq:hep4}) is 
\[
  \bigcup_{\scriptsize \begin{array}{c}
        J \subset \{1, \ldots,n-k\} \\
        1 \in J, |J| = s-k +1
   \end{array}} \left(\partial \mbox{conv} \left(e_j, j \in J\right) \times [0,1]^{QI}\right) \cup 
   \left(\mbox{conv} \left(e_j, j \in J\right) \times \{0,1\}^{QI}\right).
\]
The latter set is connected if $s-k+QI>1$. Consequently, the number of connected components of $M^{a}$ remains unchanged, and it holds:
\[
   q_a = q_{a-\varepsilon}.
\]
\end{itemize}

Now, we proceed with the {\it global argument}. Assumption \ref{ass:lb} implies that there exists $c \in \R$ such that $M^c$ is empty, thus, $q_c =0$. Additionally, there exists $d \in \R$ such that $M^d$ is connected and contains all M-stationary points, thus, $q_d=1$. Due to Assumption \ref{ass:prop}, $M_{c}^{d}$ is compact, moreover, it contains all M-stationary points. Since nondegenerate M-stationary points are in particular isolated, we conclude that there must be finitely many of them. 
Let us now increase the level $a$ from $c$ to $d$ and describe how the number $q_a$ of connected components of the lower level sets $M^a$ changes.
It follows from the local argument that $r$ new connected components are created, where $r$ is the number of local minimizers for SCNO. Let $q_I$ and $q_{II}$ denote the actual number of disappearing connected components if passing the levels corresponding to M-stationary points of types (I) and (II), respectively.
The local argument provides that at most $r_I$ and $(n-s) r_{II}$ connected components might disappear while doing so, i.\,e.
\[
  q_I \leq r_I, \quad q_{II} \leq (n-s) r_{II}.
\]
Altogether, we have:
\[
  r - r_I  - (n-s)r_{II} \leq r - q_I  - q_{II} = q_d-q_c.
\]
By recalling that $q_d=1$ and $q_c=0$, we get Morse relation (\ref{eq:mr}). \qed

We illustrate Theorem \ref{thm:mrel} by discussing the same SCNO as in Example \ref{ex:s-deg}.

\begin{example}[Saddle point]
\label{ex:s-deg1}
We consider the following SCNO with $n=2$ and $s=1$:
\[
   \min_{x_1,x_2}\,\, \left(x_1-1\right)^2 + \left(x_2-1\right)^2 \quad \mbox{s.\,t.} \quad 
   \left\|\left(x_1, x_2\right)\right\|_0 \leq 1.
\]
As we have seen in Example \ref{ex:s-deg}, both M-stationary points $(1,0)$ and $(0,1)$ are nondegenerate minimizers. Thus, we have $r=2$. Morse relation (\ref{eq:mr}) from Theorem \ref{thm:mrel} 
provides:
\[
  r_I + r_{II} \geq 1.
\]
Hence, there should exist an additional M-stationary point with M-index one. 
In fact, $(0,0)$ is this nondegenerate M-stationary point of type (II), cf. Example \ref{ex:s-deg}. Note that, due to $r_I=0$ and $r_{II}=1$, Morse relation (\ref{eq:mr}) holds with equality here.\qed
\end{example}

\section*{Acknowledgment}
The authors would like to thank Hubertus Th. Jongen for fruitful discussions.

\bibliographystyle{apalike}
\bibliography{lit.bib}

\begin{thebibliography}{}

\bibitem[Beck and Eldar, 2013]{beck:2013}
Beck, A. and Eldar, Y.~C. (2013).
\newblock Sparsity constrained nonlinear optimization: Optimality conditions
  and algorithms.
\newblock {\em SIAM Journal on Optimization}, 24:1480--1509.

\bibitem[Beck and Hallak, 2016]{beck:2016}
Beck, A. and Hallak, N. (2016).
\newblock On the minimization over sparse symmetric sets: Projections,
  optimality conditions, and algorithms.
\newblock {\em Mathematics of Operations Research}, 41:196--223.

\bibitem[Beck and Hallak, 2018]{beck:2018}
Beck, A. and Hallak, N. (2018).
\newblock Proximal mapping for symmetric penalty and sparsity.
\newblock {\em SIAM Journal on Optimization}, 28:496--527.

\bibitem[Bucher and Schwartz, 2018]{bucher:2018}
Bucher, M. and Schwartz, A. (2018).
\newblock Second-order optimality conditions and improved convergence results
  for regularization methods for cardinality-constrained optimization problems.
\newblock {\em Journal of Optimization Theory and Applications}, 178:383--410.

\bibitem[Burdakov et~al., 2016]{burdakov:2016}
Burdakov, O., Kanzow, C., and Schwartz, A. (2016).
\newblock Mathematical programs with cardinality constraints: reformulation by
  complementarity-type conditions and a regularization method.
\newblock {\em SIAM Journal on Optimization}, 26:397--425.

\bibitem[{\v{C}}ervinka et~al., 2016]{cervinka:2016}
{\v{C}}ervinka, M., Kanzow, C., and Schwartz, A. (2016).
\newblock Constraint qualifications and optimality conditions for optimization
  problems with cardinality constraints.
\newblock {\em Mathematical Programming}, 160:353--377.

\bibitem[Donoho, 2006]{donoho:2006}
Donoho, D.~L. (2006).
\newblock Compressed sensing.
\newblock {\em IEEE Transactions on Information Theory}, 52:1289--1306.

\bibitem[Dorsch et~al., 2012]{shikhman:2012a}
Dorsch, D., Shikhman, V., and Stein, O. (2012).
\newblock Mathematical programs with vanishing constraints: Critical point
  theory.
\newblock {\em Journal of Global Optimization}, 52:591--605.

\bibitem[Goerss and Jardine, 2009]{goerss:2009}
Goerss, P.~G. and Jardine, J.~F. (2009).
\newblock {\em Simplicial Homotopy Theory}.
\newblock Birkhäuser, Basel.

\bibitem[Goresky and MacPherson, 1988]{goresky:1988}
Goresky, M. and MacPherson, R. (1988).
\newblock {\em Stratified Morse Theory}.
\newblock Springer, New York.

\bibitem[Hirsch, 1976]{hirsch:1976}
Hirsch, M.~W. (1976).
\newblock {\em Differential Topology}.
\newblock Springer, Berlin-Heidelberg-New York.

\bibitem[Jongen, 1977]{jongen:1977}
Jongen, H.~T. (1977).
\newblock {\em On non-convex optimization}.
\newblock Dissertation, University of Twente, The Netherlands.

\bibitem[Jongen et~al., 2000]{Jongen:2000}
Jongen, H.~T., Jonker, P., and Twilt, F. (2000).
\newblock {\em Nonlinear Optimization in Finite Dimensions}.
\newblock Kluwer Academic Publishers, Dordrecht.

\bibitem[Jongen et~al., 2009]{shikhman:2009}
Jongen, H.~T., Shikhman, V., and Ruckmann, J.-J. (2009).
\newblock {MPCC}: Critical point theory.
\newblock {\em SIAM Journal on Optimization}, 20:473--484.

\bibitem[Milnor, 1963]{milnor:1963}
Milnor, J. (1963).
\newblock {\em Morse theory}.
\newblock Princeton University Press, Princeton, NJ.

\bibitem[Shechtman et~al., 2011]{shechtman:2011}
Shechtman, Y., Eldar, Y.~C., Szameit, A., and Segev, M. (2011).
\newblock Sparsity-based sub-wavelength imaging with partially spatially
  incoherent light via quadratic compressed sensing.
\newblock {\em Optics Express}, 19:14807--14822.

\bibitem[Shikhman, 2012]{Shikhman:2012}
Shikhman, V. (2012).
\newblock {\em Topological Aspects of Nonsmooth Optimization}.
\newblock Springer, New York.

\bibitem[Spanier, 1966]{spanier:1966}
Spanier, E.~H. (1966).
\newblock {\em Algebraic Topology}.
\newblock McGraw-Hill Book Company, New York.

\bibitem[Tibshirani, 1996]{tibshirani:1996}
Tibshirani, R. (1996).
\newblock Regression shrinkage and selection via the lasso.
\newblock {\em Journal of the Royal Statistical Society. Series B},
  58:267--288.

\end{thebibliography}

\end{document}